\documentclass{birkmult}
\usepackage[cp1251]{inputenc}
\usepackage[russian]{babel}
\usepackage{afterpage,float,amsmath,amscd}

 \newtheorem{thm}{Теорема}
 
 \newtheorem{lem}{Лемма}
 
 \theoremstyle{definition}
 
 \theoremstyle{remark}
 \newtheorem{rem}{Замечание}

\renewcommand{\labelenumi}{(\alph{enumi})}
\renewcommand{\labelenumii}{(\alph{enumii})}

\title[Orthoscalar representations of quivers]{Orthoscalar representations of quivers on the\\ category of Hilberts spaces. II.}
\author{\framebox{Roiter A.V.}\,}
\address{Институт математики НАН Украины}

\author{Kruglyak S.A.}
\address{Институт математики НАН Украины}
\email{krug@ehl.kiev.ua}
\author{Nazarova L.A.}
\address{Институт математики НАН Украины}
\email{nazarovala@yahoo.co.uk}

\begin{document}
\renewcommand{\abstractname}{Abstract}

\begin{abstract}
    As it is known, finitely presented quivers correspond to Dynkin
    graphs (Gabriel, 1972) and tame quivers -- to extended Dynkin
    graphs (Donovan and Freislich, Nazarova, 1973). In the article
    "Locally scalar reresentations of graphs in the category of
    Hilbers spaces" (Func. Anal. and Appl., 2005) authors showed the
    way to tranfer these results to Hilbert spaces, constructed
    Coxeter functors and proved an analogue of Gabriel theorem fol locally scalar
    (orthoscalar in the sequel) representations (up to the unitary
    equivalence).

    The category of orthoscalar representations of a quiver can be
    considered as a subcategory in the category of all
    representations (over a field $\mathbb C$). In the present paper
    we study the connection between indecomposable orthoscalar
    representations in the subcategory and in the category of all
    representations.

    For the quivers, corresponded to extended Dynkin graphs,
    orthoscalar representations which cannot be obtained from the
    simplest by Coxeter functors (regular representations) are
    classified.
\end{abstract}

\maketitle

\sloppy

\section*{Введение}

Более 30-ти лет назад (см. \cite{Roi75}) в работах И.М.~Гельфанда,
В.А.~Пономарева, П.~Габриеля и авторов было показано, что ряд
проблем линейной алгебры (возникшие в теории представлений алгебр,
теории групп и модулей Хариш---Чандра) допускают содержательное
изучение как на наивном языке приведения наборов матриц теми или
иными допустимыми преобразованиями, так и в категорно-функторных
терминах. Отметим представления колчанов и частично упорядоченных
множеств. Было доказано, в частности, что конечно-представимые
(ручные) колчаны соответствуют графам (расширенным графам) Дынкина.

В \cite{KruRoi05} был указан путь перенесения этих результатов на
представления колчанов в гильбертовых пространствах, построены
функторы Кокстера и доказан аналог теоремы Габриеля \cite{Gabriel72}
для локально-скалярных представлений колчанов (в дальнейшем мы
заменим термин "локально-скалярные представления"\ \cite{KruRoi05}
на термин "ортоскалярные представления"\ \!\!\!\!, считая его более
удачным). В \cite{KruRoi05} доказано, что у колчанов,
соответствующих расширенным графам Дынкина, размерности неразложимых
ортоскалярных представлений не ограничены в совокупности (как
следует из \cite{KruRabSam02}-\cite{AlbOstSam06} расширенные графы
Дынкина среди не конечно представимых и только они не имеют
бесконечномерных неразложимых ортоскалярных представлений). Мы даём
описание неразложимых ортоскалярных представлений для расширенных
графов Дынкина (см. теоремы 1-3), указывая не только на их сходство,
но и на отличие от неразложимых представлений таких графов в
линейных пространствах (см. замечание 2).


\section{О неразложимости в категории представлений колчана и её
подкатегории ортоскалярных представлений}

Напомним некоторые определения и факты (см. \cite{KruNazRoi06}) об
ортоскалярных представлениях колчанов. Колчан $Q$ с множеством
вершин $Q_v$, $|Q_v|=N$, и множеством стрелок $Q_a$ называется
\emph{разделённым}, если $Q_v=\overset{\circ}{Q}\bigsqcup
\overset{\bullet}{Q}$, и для любой $\alpha \in Q_a$ её начало
$t_\alpha \in \overset{\circ}{Q}$ и конец $h_\alpha \in
\overset{\bullet}{Q}$. Колчан $Q$ \emph{однократный}, если при
$\alpha \neq \beta$ либо $t_\alpha \neq t_\beta$, либо $h_\alpha
\neq h_\beta$. Вершины из $\overset{\circ}{Q}$ будем называть
\emph{четными}, из $\overset{\bullet}{Q}$ --- \emph{нечётными}.

Пусть $m=\big|\overset{\bullet}{Q}\big|$,
$n=\big|\overset{\circ}{Q}\big|$,
$\overset{\bullet}{Q}=\{i_1,i_2,\ldots,i_m\}$,
$\overset{\circ}{Q}=\{j_1,i_2,\ldots,j_n\}$. Представление $T$
колчана $Q$ ставит в соответствие вершине $i\in Q_v$ конечномерное
линейное пространство $T(i)$, а стрелке $\alpha: j \rightarrow i, \
\alpha \in Q_a$ линейное отображение $T_{ij}: T(j) \rightarrow
T(i)$.

Представление $T$ однократного разделённого колчана при
фиксированных базисах пространств $T(i)$, $i\in Q_v$, можно
ассоциировать с матрицей, разделённой на $m$ горизонтальных и $n$
вертикальных полос, т.е. блочной матрицей
$$
    T=\big[ T_{i_l,j_k} \big]_{k=\overline{1,n},\ l=\overline{1,m}}.
$$
При этом будем считать, что $T_{i_l,j_k}=0$, если не существует
$\alpha \in Q_a$ такой, что  $t_\alpha=j_k, h_\alpha=i_l$.

Пусть $\textrm{Rep}Q$ --- категория представлений колчана $Q$,
объекты которой есть представления, а морфизм представления $T$ в
представление $\widetilde T$ определяется как семейство линейных
отображений $C=\{C_i\}_{i\in Q_v}$,\linebreak $C_i: T(i) \rightarrow
\widetilde T(i)$, таких, что для каждой $\alpha \in Q_a$ c
$t_\alpha=j$, $h_\alpha = i$ диаграмма
\begin{equation} \label{CD}
  \begin{CD}
    T(j) @>T(\alpha)>> T(i)\\
    @VVC_jV    @VVC_iV  \\
    \widetilde T(j) @>\widetilde T(\alpha)>> \widetilde T(i)
  \end{CD}
\end{equation}
коммутативна, т.е. $C_iT_{ij}=\widetilde T_{ij}C_j$.

Пусть представления $T, \widetilde T$ заданы матрицами
$$
    T=\big[ T_{i_l,j_k} \big]_{k=\overline{1,n},\ l=\overline{1,m}}
    \quad
    \mbox{и} \quad
    \widetilde T=\big[\widetilde T_{i_l,j_k} \big]_{k=\overline{1,n},\ l=\overline{1,m}}
$$
Введём матрицы $A=\textrm{diag}\{C_{i_1},C_{i_2},\ldots,C_{i_m}\}$,
$B=\textrm{diag}\{C_{j_1},C_{j_2},\ldots,C_{j_n}\}$. Тогда из
коммутативности диаграмм (\ref{CD}) следует
\begin{equation}    \label{CD2}
    AT=\widetilde{T}B.
\end{equation}
Будем в дальнейшем говорить, что $C=(A,B).$

Пусть $\mathcal H$ --- категория унитарных (конечномерных
гильбертовых) пространств. Обозначим через $\textrm{Rep}(Q,\mathcal
H)$ подкатегорию в $\textrm{Rep}Q$, объекты которой есть
представления $T$, для которых $T(i)$ --- унитарные пространства
$(i\in Q_v)$, а морфизмы $C:T\rightarrow \widetilde T$ --- те из
морфизмов в $\textrm{Rep}Q$, для которых, кроме (\ref{CD}),
коммутативными будут и диаграммы
\begin{equation} \label{CD3}
  \begin{CD}
    T(j) @<T(\alpha)^*<< T(i)\\
    @VVC_jV    @VVC_iV  \\
    \widetilde T(j) @<\widetilde T(\alpha)^*<< \widetilde T(i)
  \end{CD}
\end{equation}
т.е. будут выполняться равенства
\begin{equation} \label{CD4}
    AT=\widetilde{T}B, \quad   BT^*=\widetilde{T}^*A
\end{equation}

Представления $T,\ \widetilde T$ из $\textrm{Rep}Q$ (соотв., из
$\textrm{Rep}(Q,\mathcal H)$) эквивалентны в $\textrm{Rep}Q$
(соотв., в $\textrm{Rep}(Q,\mathcal H)$), если найдётся обратимый
морфизм $C:T\rightarrow \widetilde T$. Можно показать, что $T$ и
$\widetilde T$ эквивалентны в $\textrm{Rep}(Q,\mathcal H)$ тогда и
только тогда, когда они унитарно эквивалентны (см., например,
\cite{Roi79}), т.е. обратимый морфизм можно выбрать состоящим из
унитарных матриц $C_i$.

Обозначим
$\overrightarrow{T_i}=\big[T_{i,j_1};T_{i,j_2};\ldots;T_{i,j_n}\big]$,
\quad $T_j^\downarrow=\left[\begin{array}{c}
                        T_{i_1,j} \\
                        T_{i_2,j} \\
                        \vdots \\
                        T_{i_m,j}
                      \end{array}\right].$

Представление $T$ разделенного однократного колчана $Q$ из категории
$\textrm{Rep}(Q,\mathcal H)$ назовём \emph{ортоскалярным}, если
каждому $i\in Q_v$ сопоставлено вещественное неотрицательное число
$\chi_i$, и выполняются следующие условия:
\begin{align} \label{OS1}
    \overrightarrow{T_i}\cdot \overrightarrow{T_i}^*&=\chi_i I_i
    \quad
    \mbox{при} \ \ i \in \overset{\bullet}{Q}, \\ \label{OS}
    T_j^{\downarrow *}\cdot T_j^{\downarrow}&=\chi_j I_j \quad
    \mbox{при} \ \ j \in \overset{\circ}{Q},
\end{align}
здесь $I_i$ --- матрица единичного оператора в $T(i)$. Ортоскалярной
будем называть и матрицу представления $T$.

В определении представлений колчана $Q$ в $\mathcal H$ можно было бы
отказатся от конечномерности пространств $T(i)$, расматривая и
бесконечномерные представления.

Будем говорить, что $Q$ (ортоскалярно) \emph{конечнопредставим} в
$\mathcal H$, если все его ортоскалярные представления распадаются в
прямую сумму (конечную либо бесконечную) конечномерных
представлений, размерности неразложимых представлений ограничены в
совокупности и в каждой размерности число неразложимых представлений
с данным характером конечно.

Определим категорию $\textrm{Rep}_{os}(Q,\mathcal H)$ как полную
подкатегорию в $\textrm{Rep}(Q,\mathcal H)$, объекты которой есть
ортоскалярные представления колчана $Q$.

Обозначим через $\mathbb R^{G_v}$ линейное вещественное протранство
из наборов $x=(x_i)$ действительных чиел $x_i\ (i\in G_v)$, элементы
$x$ из $\mathbb R^{G_v}$ будем называть $G$-векторами.

Ортоскалярному представлению $T$ разделённого однократного колчана
$Q$ сопоставим два $N$-мерных $G$-вектора ($N=m+n$):
\emph{размерность} $d=\{d(j)\}_{j\in Q_v}$ представления $T$, где
$d(j)=\dim T(j)$, и \emph{характер} $\chi=\{\chi(j)\}_{j\in Q_v}, \
\chi(j)=\chi_j$ определены выше (см. (\ref{OS1}), (\ref{OS})). Два
ортоскалярных представления $T$ и $\widetilde T$ эквивалентны, если
существуют такие унитарные матрицы
$U=\textrm{diag}\{U_{i_1},U_{i_2},\ldots,U_{i_m}\}$ и
$V=\textrm{diag}\{V_{j_1},V_{j_2},\ldots,V_{j_n}\}$, что
\begin{align*} \label{OS2}
    UT&=\widetilde TV, \quad \mbox{или} \\
    \widetilde{T}_{i_l,j_k}&=U_{i_l}T_{i_l,j_k}V_{j_k}^*
\end{align*}

Представление $T$ будем называть \emph{шуровским} (\emph{brick}) в
категории $\textrm{Rep}Q$ (соответственно $\textrm{Rep}(Q,\mathcal
H)$, $\textrm{Rep}_{os}(Q,\mathcal H)$), если его кольцо
эндоморфизмов в этой категории одномерно (изоморфно $\mathbb C$).
Очевидно, если $T$ --- шуровское представление, то $T$ неразложимо
(в соответствующей категории). Если $T$ --- неразложимое в категории
$\textrm{Rep}(Q,\mathcal H)$, то оно в ней шуровское. Действительно
алгебра $\mathfrak{A}=\textrm{End}T$ есть конечномерная $*$-алгебра.
Если $C=(A,B) \in \textrm{End}T$ то $C^*=(A^*,B^*)$. Если $C \in
\textrm{Rad} \mathfrak{A}$, то $CC^*=(AA^*,BB^*) \in \textrm{Rad}
\mathfrak{A}$ и $CC^*$ --- нильпотентный элемент, поэтому $AA^*,
BB^*$ нильпотентные и положительные операторы. Значит, $C=(0,0)$, и
алгебра $\mathfrak{A}$ полупростая. С другой стороны, алгебра
$\mathfrak{A}$
--- локальна, как алгебра эндоморфизмов неразложимого представления.
Значит, $\mathfrak{A}\simeq\mathbb C$.

Представление $T$ колчана $Q$ называется \emph{точным} если
$T(i)\neq 0$ при всех $i \in Q_v$. \emph{Носителем} представления
$T$ называется множество $Q_v^T=\{i\in Q_v \ |\  T(i)\neq 0\}$.
Характер представления определён однозначно на носителе $Q_v^T$
представления (и неоднозначно вне носителя). Если $Q_v^T=Q_v$, то
характер ортоскалярного представления определён однозначно и
обозначается $\chi_T$, в общем случае обозначим через $\{\chi_T\}$
множество всех характеров представления $T$. Ясно, что если $T$ и
$\widetilde T$ унитарно эквивалентны, то
$\{\chi_T\}=\{\chi_{\widetilde T}\}.$

Категория $\textrm{Rep}_{os}(Q,\mathcal H)$ есть (неполная)
подкатегория категории $\textrm{Rep}(Q)$. Установим связь между
неразложимыми объектами этих категорий. Основной вопрос: останется
ли неразложимое ортоскалярное представление неразложимым в категории
$\textrm{Rep}(Q)$, и поскольку, как мы отмечали, неразложимые в
$\textrm{Rep}_{os}(Q,\mathcal H)$ представления есть шуровские,
останутся ли они шуровскими в категории $\textrm{Rep}(Q)$?

Второй пункт нижеприведенного утверждения доказан в
\cite{KruNazRoi06}.

\begin{thm} \label{T1}
    Пусть $Q$ разделённый однократный колчан.
\begin{enumerate}
         \item Если $T$, $\widetilde T \in \textrm{Rep}_{os}(Q,\mathcal H)$ ---
представления с одинаковым характером и $T$ эквивалентно $\widetilde
T$ в $\textrm{Rep}(Q)$, то $T$ эквивалентно $\widetilde T$ в
$\textrm{Rep}_{os}(Q,\mathcal H)$;
         \item Если $T$ --- неразложимое представление  в $\textrm{Rep}_{os}(Q,\mathcal
         H)$ то $T$ --- неразложимое, более того, шуровское и в $\textrm{Rep}(Q)$.
       \end{enumerate}
\end{thm}

Доказательство теоремы \ref{T1} опирается на следующие
вспомогательные утверждения.

\begin{lem}
  Пусть $C=(A,B)$ --- морфизм представления $T$ в $\widetilde T$ в категории  $\textrm{Rep}(Q)$, т.е. выполняется
  равенство
$$
    AT=\widetilde{T}B,
$$
и $A, B$ есть унитарные отображения. Тогда $C=(A,B)$ есть морфизм
представления $T$ в представление $\widetilde T$ и в категории
$\textrm{Rep}_{os}(Q,\mathcal H)$ т.е. выполняются и равенства
$$BT^*=\widetilde{T}^*A.$$
\end{lem}

\begin{proof}
    Из (\ref{CD2}) следует $T^*A^*=B^*\widetilde T^*$ или, учитывая унитарность $A$ и $B$,
    имеем $T^*A^{-1}=B^{-1}\widetilde T^*$. Поэтому
    $$BT^*=\widetilde{T}^*A.$$
\end{proof}

\begin{lem}
    Пусть $C=(A,B)$ морфизм представления $T$ в себя (эндоморфизм
    представления $T$) в категории $\textrm{Rep}(Q)$, т.е.
    \begin{equation} \label{CD7}
        AT=TB
    \end{equation}
    и $A,\ B$ --- самосопряжённые операторы. Тогда $C=(A,B)$ есть
    эндоморфизм представления $T$ и в категории $\textrm{Rep}(Q, \mathcal
    H)$, т.е. и
    \begin{equation} \label{CD8}
        AT^*=T^*B.
    \end{equation}
\end{lem}

\begin{proof}
Действительно, (\ref{CD8}) получается из (\ref{CD7}) операцией
сопряжения.
\end{proof}

\begin{lem} \label{l3} (\cite{KruNazRoi06}) Пусть
$Z=\big[z_{ij}\big]_{i=1,m,\ j=1,n}$, $W=\big[w_{ij}\big]_{i=1,m,\
j=1,n}$ --- матрицы над полем $\mathbb C$, имеющие одинаковые
положительные длины ($|\overrightarrow{x}|$, $|y^\downarrow|$)
соответствующих срок и соответствующих столбцов. Пусть
$A=diag\{a_1,a_2,\ldots,a_m\}$, $B=diag\{b_1,b_2,\ldots,b_n\}$
--- матрицы над $\mathbb R$, $a_i>0,\ b_j>0$ при $i=1,m$, $j=1,n$, и
пусть $AZ=WB$. Тогда $Z=W$.
\end{lem}

\begin{proof} Обозначим через $K$ число ненулевых элементов в каждой из
матриц $Z, W$. Доказательство леммы \ref{l3} проводится
(\cite{KruNazRoi06}) индукцией по тройкам чисел $(m,n,K);$ считая,
что $(m_1,n_1,K_1)<(m_2,n_2,K_2)$, если $m_1\leq m_2$, $n_1\leq n_2$
и хотя бы одно неравенство строгое, либо если $m_1=m_2$, $n_1=n_2$,
но $K_1<K_2$. Базу индукции получаем при $m=1$, либо $n=1$, либо
$K=\max(m,n)$.
\end{proof}

\emph{Доказательство теоремы \ref{T1}.}

(a) Пусть $C=(A,B)$ осуществляет эквивалентность представлений $T$,
$\widetilde T$ из $\textrm{Rep}_{os}(Q,\mathcal H)$ в категории
$\textrm{Rep}(Q)$, т.е. $AT=\widetilde TB$ ($A$ и $B$ обратимые
матрицы). Пусть $A=XU$, $B=VY$ --- полярные разложения матриц $A$ и
$B$, где $U, V$ ---  унитарны, $X, Y$ --- положительные невырожденые
матрицы (при этом можно считать, что $U,V,X,Y$ имеют такую же
блочно-диагональную структуру, как $A$ и $B$). Пусть, кроме того,
$$
    X=U_1^*\widetilde X U_1, \quad  Y=V_1\widetilde Y V_1^*,
$$
где $U_1, V_1$ - унитарные матрицы, $\widetilde X, \widetilde Y$ ---
диагональные матрицы с положительными числами на диагонали. Тогда
\begin{equation} \label{CD9}
    U_1^*\widetilde X U_1 U T= \widetilde T V V_1 \widetilde Y V_1^*
\end{equation}
или
$$
    \widetilde X(U_1 U T V_1)= (U_1\widetilde T V V_1) \widetilde Y.
$$
Так как длины соответствующих строк и столбцов матриц $U_1UTV_1$ и
$U_1\widetilde T V V_1$ равны в силу предположений, то по лемме 3
$$
    U_1 U T V_1= U_1\widetilde T V V_1,
$$
а тогда, после сокращений,
\begin{equation} \label{CD10}
    UT=\widetilde T V,
\end{equation}
и представления $T$ и $\widetilde T$ эквивалентны в
$\textrm{Rep}_{os}(Q,\mathcal H)$ по лемме 1.

(b) Пусть $T$-шуровское представление в
$\textrm{Rep}_{os}(Q,\mathcal H)$ и $C=(A,B)\in \textrm{End} T$ в
$\textrm{Rep}(Q)$. Можно считать, что $A$ и $B$ невырождены
(прибавляя в случае необходимости к $A$ и $B$ подходящее кратное
единичного оператора). Тогда по (\ref{CD10}) при $T=\widetilde T$ из
того, что $(U,V)$ есть эндоморфизм представления $T$ в
$\textrm{Rep}_{os}(Q,\mathcal H)$, следует что $U$ и $V$ кратны
единичному оператору (с одним и тем же скаляром в качестве
сомножителя). Сокращая равенство (\ref{CD9}) на этот скаляр, мы
получим
$$
    XT=TY,
$$
где $X$ и $Y$ --- самосопряжённые операторы, а тогда по лемме 2
из-за шуровости $T$ в $\textrm{Rep}_{os}(Q,\mathcal H)$ операторы
$X$ и $Y$ скалярны (с одним и тем же скаляром); следовательно,
такими будут и операторы $A=XU$ и $B=VY$. Значит $T$ ---
неразложимое шуровское представление в $\textrm{Rep}(Q)$.

\section{Размерности неразложимых ортоскалярных представлений}

С колчаном $Q$ связана \emph{форма Титса} $q(x)$ на $\mathbb
R^{Q_v}$: если $x\in \mathbb R^{Q_v}$, то
$$
    q(x)=\sum\limits_{i\in Q_v} x_i^2-\sum\limits_{\alpha \in Q_a}
    x_{t_\alpha} x_{h_\alpha}.
$$
Из работы \cite{Kac82} следует, что размерности неразложимых (в
$\textrm{Rep}(Q)$) представлений колчана $Q$ совпадают с
положительными \emph{корнями}  соответствующего графа $G=G(Q)$,
причём, для графов Дынкина и расширенных графов Дынкина такие корни
совпадают в точности с решениями уравнений $q(x)=1$ и $q(x)=0, x\in
\mathbb Z_{+}^{Q_v}$ (здесь $\mathbb Z_{+}^{Q_v}=\{x\in \mathbb
Z^{Q_v} | x\neq0, x_i\geq0\}$). Корни $x$ при $q(x)=1$ называются
\emph{действительными} а при $q(x)=0$ --- \emph{мнимыми}. Мнимые
корни кратны минимальному мнимому положительному корню
$\delta=\delta_{G}$.

Фиксируем нумерацию вершин в
$\overset{\bullet}{Q}=\{i_1,i_2,\ldots,i_m\}$, и
$\overset{\circ}{Q}=\{j_1,i_2,\ldots,j_n\}$ (для определенности
будем считать что узловая точка, если она одна, лежит в
$\overset{\bullet}{Q}$; будем вершины из $\overset{\circ}{Q}$
обозначать также как $i_{m+k}=j_k, \ k\in \overline{i,n}$). Пусть
$x\in \mathbb R^{Q_v}$, $x_k=x(i_k)$, при $k\in \overline{1,m+n}$,
$c$---преобразование Кокстера на $\mathbb R^{Q_v}$,
$c=\sigma_{i_{n+m}}\cdots\sigma_{i_2}\sigma_{i_1}$,
$(\sigma_{i_k}(x))_k=-x_k+\sum\limits_{l,i_l-i_k}x_l$,
$(\sigma_{i_k}(x))_l=x_l$ при $l\neq k$. Ясно, что
$\sigma_i^2=\textrm{id}$ при $i\in \overline{1,n+m}$. Поэтому
$c^{-1}=\sigma_{i_1}\sigma_{i_2}\cdots\sigma_{i_{n+m}}$. Будем также
пользоваться обозначениями
$\overset{\bullet}{c}=\sigma_{i_{n}}\cdots\sigma_{i_2}\sigma_{i_1}$,
$\overset{\circ}{c}=\sigma_{i_{m+n}}\cdots\sigma_{i_{m+2}}\sigma_{i_{m+1}}$
и называть эти преобразования отражениями Кокстера
(${\overset{\bullet}{c}{}^2}=\textrm{id}$, $\overset{\circ}{c}
{}^2=\textrm{id}$). Обозначим
$\overset{\bullet}{c}_k=\underbrace{\ldots
\overset{\bullet}{c}\overset{\circ}{c}\overset{\bullet}{c}}_{k \
\mbox{\Small раз}}$, $\overset{\circ}{c}_k=\underbrace{\ldots
\overset{\circ}{c}\overset{\bullet}{c}\overset{\circ}{c}}_{k \
\mbox{\Small раз}}$, $k \in \mathbb N$. Вектор $x \in \mathbb
R_{+}^{G_v}$ \emph{регулярен}, если $c^t(x) \in \mathbb R_{+}^{G_v}$
при любом $t \in \mathbb Z$  и \emph{сингулярен} в противном случае
(терминология восходит к \cite{GelPon70}).

Представление $T$ колчана $Q$ \emph{сингулярно}, если $T$
неразложимо, конечномерно и его размерность $d$ --- сингулярный
вектор; $T$ \emph{регулярно}, если $T$ неразложимо, конечномерно и
не сингулярно.

Пусть
$\delta_G=(\delta_1,\ldots,\delta_m,\delta_{m+1},\ldots,\delta_{m+n})$.
Построим линейную форму
$$
    L_{G}(x)=\sum\limits_{i_k\in \overset{\bullet}{G}} \delta_k x_k-
    \sum\limits_{i_{m+k}\in \overset{\circ}{G}} \delta_{m+k}
    x_{m+k}, \quad x\in \mathbb R^{G_v}.
$$

Справедливо утверждение (см., например, \cite{RedRoi04}): пусть $G$
--- расширенный граф Дынкина; для того чтобы корень $x\in \mathbb
R^{G_v}$ был сингулярен, необходимо и достаточно, чтобы $L_G(x)\neq
0$.

Отметим, что у графа Дынкина все корни дейтвительны и сингулярны, у
расширенного графа Дынкина все мнимые корни регулярные, а вот
действительные корни разбиваются на два сорта (действительные
сингулярные и действительные регулярные). Для колачана $Q$, граф
которого $G$ есть расширенный граф Дынкина, справедливо утверждение:

\begin{lem} (cм. \cite{CrawBoe}) Если $T$ есть неразложимое в
$\textrm{Rep}(Q)$ регулярное представление колчана $Q$, граф $G$
которого есть расширенный граф Дынкина, то $T$ --- шуровское
представление тогда и только тогда когда $\dim T=d\leq\delta_G$.
\end{lem}

Докажем следующее утверждение:

\begin{lem}
Пусть $G$ --- расширенный граф Дынкина, $d\in \mathbb Z_{+}^G$
корень графа $G$, удовлетворяющий условию $d<\sigma_G$. Тогда либо
$d$ --- неточный корень (одна из координат нулевая), либо получается
из некоторого неточного корня $\tilde d$ применением отражений
Кокстера ${\overset{\bullet}{c}}$ и ${\overset{\circ}{c}}$.
\end{lem}

\begin{proof}
    Достаточно показать, что такой точный корень можно отражениями
    Кокстера преобразовать в неточный корень $\widetilde d.$

    Рассмотрим граф $\widetilde D_n$

\begin{center}
\begin{picture}(140,50)

\put(0,40){\circle*{3}} \put(20,20){\circle*{3}}
\put(50,20){\circle*{3}} \put(60,18){\circle*{1}}
\put(70,18){\circle*{1}} \put(80,18){\circle*{1}}
\put(90,20){\circle*{3}} \put(120,20){\circle*{3}}
\put(140,40){\circle*{3}} \put(140,0){\circle*{3}}
\put(0,0){\circle*{3}}

\put(0,40){\line(1,-1){20}} \put(0,0){\line(1,1){20}}
\put(20,20){\line(1,0){30}}\put(90,20){\line(1,0){30}}
\put(140,40){\line(-1,-1){20}} \put(140,0){\line(-1,1){20}}

\put(-15,40){$a_1$} \put(-15,0){$a_2$} \put(20,25){$z_1$}
\put(50,25){$z_2$} \put(90,25){$z_{n-4}$} \put(125,20){$z_{n-3}$}
\put(155,40){$b_1$} \put(155,0){$b_2$}
\end{picture}
\end{center}
\medskip


    Пусть
    $d=(d_{a_1},d_{a_2},d_{z_1},\ldots,d_{z_{n-3}},d_{b_1},d_{b_2}).$
    Как известно,
    $$
        \delta_{\widetilde{D}_n}=\begin{array}{c}
                                   1 \\
                                   1
                                 \end{array}
         2\ 2\ \ldots\ 2\ 2 \begin{array}{c}
                                   1 \\
                                   1
                                 \end{array}
    $$
    (координаты вектора располагаются так же, как соответствующие
    вершины в графе).
    Так как $d$ --- точный корень, то
    $d_{a_1}=d_{a_2}=d_{b_1}=d_{b_2}=1$. Назовём \emph{отмеченным} первый
    индекс $k$, для которого $x_{z_k}=1$ (т.е.
    $x_{z_1}=x_{z_2}=\ldots=x_{z_{k-1}}=2)$). Пусть, для
    определённости, $z_{k-1} \in \overset{\bullet}{G}$. Тогда у
    $\overset{\bullet}{c}d$ координаты
    $x_{a_1},x_{a_2},x_{z_1},\ldots,x_{z_{k-2}}$ совпадают с
    соответствующими координатами вектора $d$, а $x_{z_{k-1}}=1$,
    т.е. отмеченный индекс уменьшается на единицу. Применяя к
    $\overset{\bullet}{c}d$ преобразование
    $\overset{\circ}{c}$, продолжим этот процесс. В результате мы
    придём к вектору, у которого отмеченный индекс равен 1, а стало
    быть на следующем шаге получим вектор, у которого
    $x_{a_1}=x_{a_2}=0$, это и будет искомый неточный вектор $\tilde
    d$.

    Для графов $\widetilde E_6$, $\widetilde E_7$, $\widetilde E_8$
    результат можно получить аналогичными рассуждениями или, в
    крайнем случае, прямым перебором (положительных корней $d$ таких,
    что $d<\delta_G$, конечное число).

     $\delta_{\widetilde
    A_n}=(1,1,\ldots,1)$, так что из условия $d<\delta_{\widetilde
    A_n}$ следует, что одна из координат нулевая и $d$ есть неточный
    вектор.
\end{proof}

\section{Функторы Кокстера и ортоскалярные представления, размерности
которых есть действительные корни}

Пусть $g\in Q_v$ и $\Pi_g$ --- простейшее представление колчана $Q:
\Pi_g(g)=\mathbb C$, $\Pi_g(i)=0$ при $i\neq g,\ i\in Q_v$. Ясно,
что если $f_g$ --- характер представления $\Pi_g$, то $f_g(g)=0$;
будем предполагать, что $f_g(i)>0$ при $i\neq g$.

Обозначим через $\textrm{Rep}(Q,d,\chi)$ полную подкатегорию
неразложимых представлений из $\textrm{Rep}_{os}(Q,\mathcal H)$ с
фиксироваными размерностью $d$ и характером $\chi$. Если в
$\textrm{Rep}(Q,d,\chi)$ не входит простейшее представление, то
$\chi(i)>0$ для $d(i)\neq 0$ (в силу неразложимости представлений).

В \cite{KruRoi05} были введены \emph{функторы отражений Кокстера}
$\overset{\bullet}{F}$ и $\overset{\circ}{F}$:
\begin{equation*}
\begin{split}
    \overset{\circ}{F}: \textrm{Rep}(Q,d,\chi) \rightarrow
                        \textrm{Rep}(Q,\overset{\circ}{c}(d),\overset{\circ}{\chi}),
     \\
    \overset{\bullet}{F}: \textrm{Rep}(Q,d,\chi) \rightarrow
                        \textrm{Rep}(Q,\overset{\bullet}{c}(d),\overset{\bullet}{\chi}),
\end{split}
\end{equation*}
где

\begin{equation}\label{CD11}
\begin{split}
    \overset{\circ}{d}(i)=\overset{\circ}{c}(d)(i)=
    \left\{\begin{array}{c}
          -d(i)+\sum\limits_{j:j-i}d(j) \quad \mbox{при} \ i\in \overset{\circ}{Q} \\
          d(i) \quad \mbox{при} \ i\in \overset{\bullet}{Q}
         \end{array} \right.,\\
    \overset{\bullet}{d}(i)=\overset{\bullet}{c}(d)(i)=
    \left\{\begin{array}{c}
          -d(i)+\sum\limits_{j:j-i}d(j) \quad \mbox{при} \ i\in \overset{\bullet}{Q} \\
          d(i) \quad \mbox{при} \ i\in \overset{\circ}{Q}
         \end{array} \right.,
\end{split}
\end{equation}

\begin{equation}\label{CD12}
\begin{split}
    \overset{\circ}{\chi}(i)=
    \left\{\begin{array}{c}
          -\chi(i)+\sum\limits_{j:j-i}\chi(j) \quad \mbox{при} \ i\in \overset{\bullet}{Q_v^T} \\
          \chi(i) \quad \mbox{при} \ i\notin \overset{\bullet}{Q_v^T}
         \end{array} \right.,\\
    \overset{\bullet}{\chi}(i)=
    \left\{\begin{array}{c}
          -\chi(i)+\sum\limits_{j:j-i}\chi(j) \quad \mbox{при} \ i\in \overset{\circ}{Q_v^T} \\
          \chi(i) \quad \mbox{при} \ i\notin \overset{\circ}{Q_v^T}
         \end{array} \right.,\\
\end{split}
\end{equation}
(здесь $\overset{\bullet}{Q_v^T}=Q_v^T\cap\overset{\bullet}{Q}$ и
$\overset{\circ}{Q_v^T}=Q_v^T\cap\overset{\circ}{Q}$).

$\overset{\bullet}{F}$ и $\overset{\circ}{F}$ есть функторы
эквивалентности категорий, и можно проверить, что их двукратное
применение приводит к функтору, эквивалентному тождественному.

В дальнейшем будем пользоваться обозначениями
$\overset{\circ}{F}_k=\underbrace{\ldots
\overset{\circ}{F}\overset{\bullet}{F}\overset{\circ}{F}}_{k \
\mbox{\Small раз}}$, $\overset{\bullet}{F}_k=\underbrace{\ldots
\overset{\bullet}{F}\overset{\circ}{F}\overset{\bullet}{F}}_{k \
\mbox{\Small раз}}$.

\begin{rem}
    В \cite{KruRoi05} функторы отражений Кокстера вводятся в
    предположении, что граф $G$ колчана $Q$ не содержит циклов.
    Конструкция функторов $\overset{\bullet}{F}$,
    $\overset{\circ}{F}$ практически без изменений переносится на все
    однократные разделённые колчаны, и для них имеют место формулы
    (\ref{CD11}-\ref{CD12}).
\end{rem}

Пусть $T\in \textrm{Rep}_{os}(Q,\mathcal H)$, $Q_v^T$ --- носитель,
$d_T=d_T(i)$ --- размерность, $\{\chi_T=\chi_T(i)\}$ --- характер
представления $T$ (напомним, что вне $Q_v^T$ он определён
неоднозначно). Будем придавать $\chi_T(i)$ при $i\in Q_v \setminus
Q_v^T$ произвольные положительные значения, таким образом $\chi_T$
будет зависеть от $|Q_v|-|Q_v^T|=r$ положительных параметров.

Докажем теорему,

\begin{thm}
    Если $Q$ есть разделённый однократный колчан, граф $G$ которого
    есть граф Дынкина, либо расширенный граф Дынкина, $d=\{d_i\}_{i\in Q_v}$ его \emph{точный}
    (т.е. $d_i>0$) \emph{действительный корень}, то неразложымые
    ортоскалярные представления колчана $Q$ в размерности $d$
    зависят от $|Q_v|-1$ положительных параметров.
\end{thm}

\begin{proof}
    a)    Пусть $d$ ---  действительный сингулярный корень графа $G$ и
    $T$  --- неразложимое ортоскалярное представление размерности $d$, $d$
    получается из некоторого простейшего корня
    $d_g=(0,0,\ldots,1,0,\ldots,0)$ отражениями Кокстера
    $\overset{\circ}{c}, \overset{\bullet}{c}$ (1 соответствует
    вершине $g$). Пусть, для определённости,
    $g\in\overset{\bullet}{G}$. Тогда для некоторого $k\in \mathbb
    N$, $d=\overset{\circ}{c}_kd_g$, а
    $T=\overset{\circ}{F}_k(\Pi_g)$. $\Pi_g$ зависит от $|Q_v|-1$
    положительных параметров (положительных значений характера в
    вершинах $i\neq g$). От них же зависят скалярные
    операторы $\overrightarrow{T}\cdot\overrightarrow{T}^*$ при
    $i\in \overset{\bullet}{Q}$ и $T_j^{\downarrow*}\cdot T_j^{\downarrow}$ при $j\in
    \overset{\circ}{Q}$ в ненулевых пространствах. Соответствие
    между наборами параметров и неразложимыми ортоскалярными
    представлениями в размерности $d$ взаимнооднозначное в силу
    обратимости преобразований $\overset{\circ}{c}$,
    $\overset{\bullet}{c}$, $\overset{\circ}{\chi}$,
    $\overset{\bullet}{\chi}$ и $\overset{\circ}{F}$,
    $\overset{\bullet}{F}$.

    б) Пусть $d$ --- действительный регулярный корень графа $G$
    (для расширенного графа Дынкина) и $T$ --- неразложимое ортоскалярное представление размерности $d$.
    Если $d$ --- точный корень, то
    по лемме 5 он получается из некоторого неточного корня $\tilde
    d$ отражениями Кокстера $\overset{\circ}{c}$,
    $\overset{\bullet}{c}$. Пусть (для определённости) $d=\overset{\circ}{c}_k \tilde
    d$. Ортоскалярное представление $\widetilde T$ в размерности
    $\tilde d$ можно рассматривать как точное представление колчана,
    соответствующего графу Дынкина. Если $Q_v^{\widetilde T}$ ---
    носитель представления $\widetilde T$, то $\widetilde T$ зависит
    по вышеизложенному от $|Q_v^{\widetilde T}|-1$ положительных
    параметров. Значения характера представления $\widetilde T$ в
    точках $Q_v \setminus Q_v^{\widetilde T}$ положим равными произвольным
    положительным числам: $T=\overset{\circ}{F}_k(\widetilde T)$, представление $T$ зависит от $|Q_v \setminus Q_v^{\widetilde T}|+|Q_v^{\widetilde T}|-1=|Q_v|-1$
    положительных параметров.
\end{proof}

\section{Ортоскалярные представления, размерности которых есть
мнимые корни}

Пусть $Q$ --- разделённый однократный колчан, граф которого $G$ есть
расширенный граф Дынкина. Из теоремы 1 и леммы 4 следует, что
единственный мнимый корень, который может быть размерностью
неразложимого ортоскалярного представления колчана $Q$ есть
$\delta_G$.

Докажем теорему аналогичную теореме 2, но относящуюся к мнимым
корням графа $G$ (заметим, что в случае $\chi_T=\delta_G$ для графов
$\widetilde D_4, \widetilde E_6, \widetilde E_7$ на другом языке и
другими методами унитарная классификация неразложимых представлений
получена в \cite{OstSam99}-\cite{Ost04}).

\begin{thm} \label{T3}
    Если $Q$ есть разделённый однократный колчан, граф $G$ которого
    есть расширенный граф Дынкина, то в размерности $\delta_G$
    неразложимые ортоскалярные представления колчана $Q$ зависят от
    $|Q_v|+1$ положительных параметров.
\end{thm}

Теорему докажем фактическим описанием всех неразложимых
представлений указанных колчанов в размерностях $\delta_G$.

a) Пусть $G=\widetilde A_n$\! ---\! цикл с $n$ вершинами. Из
предположения о разде\-лен\-но\-сти и однократности колчана следует,
что $n$ чётное и $n\!\geq4$. $\delta_{\widetilde
A_n}=(1,1,\ldots,1).$ Представление $T$ в размерности
$\delta_{\widetilde A_n}$ задаётся $n$ ненулевым комплексными
числами. Переходом к унитарноэквивалентному представлению $n-1$
число можно сделать положительным ($|Q_v|-1$ параметр); допустимые
преобразования, не меняющие значения этих параметров, не меняют и
$n$--е число (зависящее от 2-х положительных параметров
--- аргумента и модуля). Таким образом, $T$ зависит от $n$
параметров, принимающих произвольные положительные значения, и
одного параметра $\phi$, принимающего значения на отрезке
$(0,2\pi]$, т.е. от $|Q_v|+1=n+1$ параметров.

\medskip

b) Пусть $G=\widetilde D_n \ (n\geq4)$:

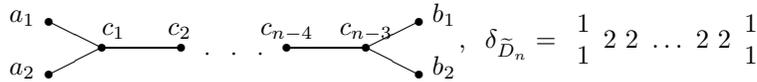
\begin{center}
\begin{picture}(240,40)
\put(0,30){\circle*{3}} \put(20,20){\circle*{3}}
\put(50,20){\circle*{3}} \put(60,18){\circle*{1}}
\put(70,18){\circle*{1}} \put(80,18){\circle*{1}}
\put(90,20){\circle*{3}} \put(120,20){\circle*{3}}
\put(140,30){\circle*{3}} \put(140,10){\circle*{3}}
\put(0,10){\circle*{3}}

\put(0,30){\line(2,-1){20}} \put(0,10){\line(2,1){20}}
\put(20,20){\line(1,0){30}}\put(90,20){\line(1,0){30}}
\put(140,30){\line(-2,-1){20}} \put(140,10){\line(-2,1){20}}

\put(-15,30){$a_1$} \put(-15,10){$a_2$} \put(20,25){$c_1$}
\put(45,25){$c_2$} \put(80,25){$c_{n-4}$} \put(110,25){$c_{n-3}$}
\put(145,30){$b_1$} \put(145,10){$b_2$}

\put(155,20){\mbox{$,$}}

\put(165,20){$   \delta_{\widetilde D_n}=\begin{array}{c}
                              1 \\
                              1
                            \end{array}
     2\ 2\ \ldots\ 2\ 2 \begin{array}{c}
                              1 \\
                              1
                            \end{array}$}

\end{picture}
\end{center}
\medskip

 Пусть $T$ --- неразложимое ортоскалярное
представление размерности $\delta_{\widetilde D_n}$ соответствующего
разделённого колчана $Q$. Для определённости, будем считать, что
$c_1 \in \overset{\bullet}{Q}$.

Для удобства введём обозначения и нумерацию матриц представления
следующим образом:
\begin{equation*}
\begin{split}
    A_1&=T_{c_1,a_1}, \  A_2=T_{c_1,a_2}, \ X_1=T_{c_1,c_2}, \
    X_2=T_{c_3,c_2},\ \ldots \\
    B_1&=T_{c_{n-3},b_1}, \  B_2=T_{c_{n-3},b_2} \quad \mbox{при}\ \
    n \ \
    \mbox{чётном и}  \\
    B_1&=T_{b_1,c_{n-3}}, \  B_2=T_{b_2,c_{n-3}} \quad \mbox{при} \ \
    n \ \
    \mbox{нечётном.}
\end{split}
\end{equation*}
Рассмотрим случай чётного $n$. При $n=4$ ортоскалярные неразложымые
представления графа $Q$ в размерности $\delta_{\widetilde D_4}$
зависят (см. \cite{KruNazRoi06}) от $6=|Q_v|+1$ положительных
параметров.

Пусть $n>4$. Допустимыми (унитарными) преобразованиями матриц
представления матрицы $X_1, X_2,\ldots,X_{n-4}$ можно
диагонализировать с неотрицательными числами на диагоналях. Пусть
$A=[A_1|A_2]$, $X_i=\left[
                            \begin{array}{cc}
                              x_i & 0 \\
                              0 & y_i \\
                            \end{array}
                          \right]$, $B=\left[
                                                  \begin{array}{c|c}
                                                    B_1 & B_2 \\
                                                  \end{array}
                                                \right]
                          $,
                          $i=\overline{1,n-4}$. Тогда из условий
 ортоскалярности представления мы легко получаем, что матрицы $AA^*$
 и $BB^*$ также диагональны. Пусть
 \begin{equation} \label{CD13}
    AA^*=\left[
                            \begin{array}{cc}
                              x^2_0 & 0 \\
                              0 & y^2_0 \\
                            \end{array}
                          \right], \quad
    BB^*=\left[
                            \begin{array}{cc}
                              x^2_{n-3} & 0 \\
                              0 & y^2_{n-3} \\
                            \end{array}
                          \right].
 \end{equation}
Тогда условия ортоскалярности выглядят так:
$$
                        \left[
                            \begin{array}{cc}
                              x^2_0 & 0 \\
                              0 & y^2_0 \\
                            \end{array}
                          \right]+
                          \left[
                            \begin{array}{cc}
                              x^2_{i+1} & 0 \\
                              0 & y^2_{i+1} \\
                            \end{array}
                          \right]=
                          \chi_{c_{i+1}}I_{c_{i+1}}, \quad i
                          \in \overline{0,n-4},
$$
и поэтому
$
    x_i^2+x_{i+1}^2=y_i^2+y_{i+1}^2,
$
а значит
\begin{equation} \label{CD14}
    y_{i+1}^2=x_{i+1}^2+(-1)^{i}(x_0^2-y_0^2), \quad i\in
    \overline{0,n-4}.
\end{equation}
Если $x_0^2=y_0^2$ (либо $x_j^2=y_j^2$ для любого фиксированого
$j$), то из (\ref{CD14}) следует, что $x_i^2=y_i^2$ для $i \in
\overline{0,n-3}$. В этом случае матрица
\begin{equation}    \label{CD15}
        \left[
      \begin{array}{c|c|c|c}
        A_1 & A_2 & B_1 & B_2
      \end{array}
    \right]
\end{equation}
задаёт ортоскалярное представление $\widehat{T}$ колчана
$\widehat{Q}$, соответствующего графу $\widetilde D_4$, и при
переходе к эквивалентному представлению колчана $Q$, не меняющему
приведенного вида матриц $X_i$, мы получаем переход к представлению
колчана $\widehat{Q}$, эквивалентному представлению $\widehat{T}$.
Такие представления колчана $\widehat{Q}$, как мы отметили выше,
зависят не более, чем от 6 параметров. Соответствующие представления
колчана $Q$ (учитывая параметры $x_1, x_2,\ldots,x_{n-4}$ и
(\ref{CD14})) зависят не более чем от $6+(n-4)=n+2=|Q_v|+1$
параметров.

Пусть $x_i\neq y_i$ при $i\in \overline{0,n-3}$ и $(U_i)_{i\in G_v}$
осуществляют эквивалентность представления $T$ c ортоскалярным
представлением $\widetilde T$, для которого $\widetilde X_i=X_i$,
$i\in \overline{0,n-4}$ ($U_i$ --- унитарные матрицы). В этом случае
легко убедится, что
$$
        U_i=\left[
              \begin{array}{c c}
                u_i & 0 \\
                0 & v_i \\
              \end{array}
            \right]
$$

Такими преобразованиями в матрице $    \left[
      \begin{array}{c|c}
        A & B
      \end{array}
    \right]=\left[
                                           \begin{array}{c|c|c|c}
                                             a_{11} & a_{12} & b_{11} &
                                             b_{12}\\
                                             a_{21} & a_{22} & b_{21} &
                                             b_{22}\\
                                           \end{array}
                                         \right]$ можно все
                                         елементы,
кроме 2-х, для определённости, $b_{21}$ и $b_{22}$, сделать
вещественными, и параметризовать следующим образом:
$$[A|B]=\left[
                                           \begin{array}{c|c|c|c}
                                             x_0\cos\phi_1 & x_0\sin\phi_1 & x_{n-3}\cos\phi_2 &
                                             x_{n-3}\sin\phi_2\\
                                             y_0\sin\phi_1 & -y_0\cos\phi_1  & y_{n-3}\sin\phi_2e^{i\theta} &
                                             -y_{n-3}\cos\phi_2e^{i\theta}\\
                                           \end{array}
                                         \right]$$
(здесь мы воспользовались тем, что
$$
    AA^*=\left[
                            \begin{array}{cc}
                              x^2_0 & 0 \\
                              0 & y^2_0 \\
                            \end{array}
                          \right], \quad \mbox{и} \quad
    BB^*=\left[
                            \begin{array}{cc}
                              x^2_{n-3} & 0 \\
                              0 & y^2_{n-3} \\
                            \end{array}
                          \right] ).
$$

Такие представления зависят от $n+2=|Q_v|+1$ положительных
параметров $x_0,\ x_1,\ \ldots,x_{n-3},\ y_0,\ \phi_1,\ \phi_2,\
\theta$.

Случай нечётного $n$ изучается аналогично с заменой матрицы
$
    \left[
      \begin{array}{c|c|c|c}
        A_1 & A_2 & B_1 & B_2
      \end{array}
    \right]
$ на матрицу
$
    \left[
      \begin{array}{c|c|c|c}
        A_1 & A_2 & B_1^* & B_2^*
      \end{array}
    \right]
$

\medskip

с) Пусть $G=\widetilde E_6$; соответствующий разделенный колчан $Q$
имеет вид

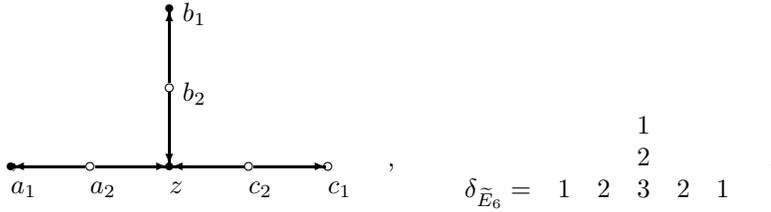
\begin{center}
\begin{picture}(160,80)

\put(0,0){\circle*{3}} \put(30,0){\circle{3}}
\put(60,0){\circle*{3}} \put(90,0){\circle{3}}
\put(120,0){\circle{3}} \put(60,30){\circle{3}}
\put(60,60){\circle*{3}}

\put(28,0){\vector(-1,0){27}} \put(32,0){\vector(1,0){27}}
\put(88,0){\vector(-1,0){27}} \put(92,0){\vector(1,0){27}}
\put(60,28){\vector(0,-1){27}} \put(60,32){\vector(0,1){27}}

\put(0,-10){$a_1$} \put(30,-10){$a_2$} \put(60,-10){$z$}
\put(90,-10){$c_2$} \put(120,-10){$c_1$} \put(65,25){$b_2$}
\put(65,55){$b_1$}, \put(140,0){\mbox{,}}

\end{picture}
\mbox{ $\begin{array}{cccccc}
                                     & &  & 1 &  &  \\
                                     & &  & 2 &  & \\
          \delta_{\widetilde E_6}=   & 1 & 2 & 3 & 2 & 1
                                   \end{array}
$\quad . \vspace*{1cm}}
\end{center}

Пусть $T$ --- неразложимое ортоскалярное представление колчана $Q$ в
размерности $\delta_{\widetilde E_6}$,
$$
    T=\begin{array}{|c|c|c|}
        \hline
        0 & 0 & C_1 \\
        \hline
        0 & B_1 & 0 \\
        \hline
        A_1 & 0 & 0 \\
        \hline
        A_2 & B_2 & C_2\\
        \hline
      \end{array}
$$
Здесь
$
    A_1=T_{a_1,a_2}; \ B_1=T_{b_1,b_2};\ C_1=T_{c_1,c_2}; \
    A_2=T_{z,a_2}; \ B_2=T_{z,b_2}; \ C_2=T_{z,c_2}.
$
Матрица $A_1$ имеет размерность $1\times2$. Унитарными
преобразованиями столбцов привёдем её к виду
$
    A_1=[0\ \ x_0], \quad x_0>0.
$
Унитарными преобразованиями строк матрицу $A_2=\left[
                                                 \begin{array}{cc}
                                                   a_{11} & a_{12} \\
                                                   a_{21} & a_{22} \\
                                                   a_{31} & a_{32} \\
                                                 \end{array}
                                               \right]$ можно
привести к виду $A_2=\left[
                                              \begin{array}{cc}
                                                   x_1 & a_{12} \\
                                                   0 & a_{22} \\
                                                   0 & a_{32} \\
                                                 \end{array}
                                               \right], \quad x_1>0 $
и из ортогональности столбцов матрицы $\begin{array}{|c|}
                                           \hline A_1 \\
                                           \hline A_2 \\
                                           \hline
                                         \end{array}$ следует
                                         $a_{12}=0$.
$$
    D=\left[
        \begin{array}{c|c|c}
          A_2 & B_2 & C_2 \\
        \end{array}
      \right]=
      \left[
        \begin{array}{cc|cc|cc}
            x_1 & 0 & b_{11} & b_{12} & c_{11} & c_{12} \\
            0 & a_{22} & b_{21} & b_{22} & c_{21} & c_{22} \\
            0 & a_{32} & b_{31} & b_{32} & c_{31} & c_{32} \\
        \end{array}
      \right].
$$
Унитарными преобразованиями столбцов (не меняя $A_2$) можно
добиться, чтобы $b_{11}=0$, $c_{11}=0$; $b_{12}\neq0$,
$c_{12}\neq0$, иначе, можно проверить, представление разложимо;
допустимыми преобразованиями можно сделать $b_{12}$ и $c_{12}$
положительными. Унитарными преобразованиями двух последних строк
матрицы $D$ можно сделать $b_{22}=0$, а тогда из ортогональности
первых двух строк получим, что $c_{22}=0$. Итак,
$$
    D=\left[
        \begin{array}{c|c|c}
          A_2 & B_2 & C_2 \\
        \end{array}
      \right]=
      \left[
        \begin{array}{cc|cc|cc}
            x_1 & 0 & 0 & b_{12} & 0 & c_{12} \\
            0 & a_{22} & b_{21} & 0 & c_{21} & 0 \\
            0 & a_{32} & b_{31} & b_{32} & c_{31} & c_{32} \\
        \end{array}
      \right].
$$

Значения характера $\chi_i$ $(i\in Q_v)$ связаны соотношением
$$
    \chi_{a_1}+\chi_{b_1}+\chi_{c_1}+3\chi_{x}=2(\chi_{a_2}+\chi_{b_2}+\chi_{c_2})
$$
(следует сумму квадратов модулей матричных элементов матрицы $T$
подсчитать двумя способами --- сначала складывая эти квадраты
модулей по строкам, а потом по столбцам). Это дает возможность
перейти к нормированному представлению со значением характера
$\chi_z=1$, поделив все матрицы представления на одно и то же число
$\sqrt{\chi_z}$ (нормировка представления не меняет минимального
числа параметров, от которых зависит представление).

Вектор $(x_1,b_{12},c_{12})$, где $x_1>0,\ b_{12}>0,\ c_{12}>0$ и
$x_1^2+b_{12}^2+c_{12}^2=1$, может быть параметризован следующим
образом
$$
    (x_1,b_{12},c_{12})=(\sin\psi_1,
    \sin\phi_1\cos\psi_1,\cos\phi_1\cos\psi_1),
$$
где $0<\phi_1<\frac{\pi}{2}$,  $0<\psi_1<\frac{\pi}{2}$. Унитарными
преобразованими столбцов элементы $a_{22}, b_{21}, c_{21}$ могут
быть сделаны неотрицательными. Так как $a_{22}^2+
b_{21}^2+c_{21}^2=1$, то вектор $(a_{22}, b_{21}, c_{21})$ может
быть параметризован следующим образом:
$$
    (a_{22},b_{21},c_{21})=(\cos\phi_2 \cos\psi_2,
    \sin\phi_2\cos\psi_2,\sin\psi_2),
$$
где $0\leq\phi_2\leq\frac{\pi}{2}$,  $0\leq\phi_2\leq\frac{\pi}{2}$,
$\phi_2\cdot\psi_2\neq 0$ (иначе представление разложимо).

Умножением 3-й строки матрицы $D$ на элемент вида $e^{i\phi}$
(унитарное преобразование строк) сделаем элемент $c_{31}$
вещественным неотрицательным
$$
    a_{32}\overline{a}_{32}+b_{31}\overline{b}_{31}+
        b_{32}\overline{b}_{32}+c_{31}^2+
        c_{32}\overline{c}_{32}=1
$$
Пусть
\begin{equation} \label{CD16}
    a_{32}\overline{a}_{32}+b_{31}\overline{b}_{31}+c_{31}^2=\sin^2\psi_3,
    \ \
    b_{32}\overline{b}_{32}+c_{31}^2=\cos^2\psi_3.
\end{equation}
Учитывая ортогональность 1-й и 3-й строки матрицы $D$ и
(\ref{CD16}), вектор $(b_{32},c_{33})$ может быть пропараметризован
так:
$$
    (b_{32},c_{33})=(-\cos\phi_1\cos\psi_3 e^{i\theta_1},\sin\psi_1\cos\psi_3
    e^{i\theta_1})
$$
Учитывая (\ref{CD16}) вектор $(a_{32}, b_{31}, c_{31})$ может быть
пропараметризован так:
$$
    (a_{32}, b_{31}, c_{31})=(\cos\phi_3\cos\psi_4\sin\phi_3 e^{i\theta_2},
                             \sin\phi_3\cos\psi_4\sin\phi_3
                             e^{i\theta_3},
                             \sin\psi_4\sin\psi_3).
$$

Таким образом, в матрице $D$, которую мы будем называть
\emph{основой представления},
\begin{equation*}
\begin{split}
A_2&=
    \left[
      \begin{array}{cc}
        \sin\psi_1 & 0\\
        0 & \cos\phi_2\cos\psi_2\\
        0 & \cos\phi_3\cos\psi_4\sin\psi_3 e^{i\theta_2} \\
      \end{array}
    \right]
\\
B_2&=
    \left[
      \begin{array}{cc}
        0 & \sin\phi_1\cos\psi_1\\
        \sin\phi_2\cos\psi_2 & 0\\
        \sin\phi_3\cos\psi_4\sin\psi_3 e^{i\theta_3} & -\cos\phi_1\cos\psi_3e^{i\theta_1} \\
      \end{array}
    \right] \\
C_2&=
    \left[
      \begin{array}{cc}
        0 & \cos\phi_1\\
        \sin\phi_2 & 0\\
        \sin\psi_4\sin\psi_3 & \sin\phi_1\cos\psi_3e^{i\theta_1} \\
      \end{array}
    \right]
\end{split}
\end{equation*}
Следовательно, основа представления, зависит от вещественных
параметров
\begin{equation} \label{CD17}
    \phi_1, \ \phi_2, \ \phi_3,\ \psi_1, \ \psi_2, \ \psi_3, \
    \psi_4, \ \theta_1, \ \theta_2, \ \theta_3.
\end{equation}
Из условий ортоскалярности cледует, что строки матрицы $D$ должны
быть ортонормированы. Из ортогональности 2-й и 3-й строки имеем
соотношение:
\begin{equation*}
\begin{split}
    &\cos\phi_2\cos\phi_3\cos\psi_2\sin\psi_3\cos\psi_4e^{i\theta_2}+
    \sin\phi_2\sin\phi_3\cos\psi_2\sin\psi_3\cos\psi_4e^{i\theta_3}+\\&+
    \sin\psi_2\sin\psi_3\sin\psi_4=0,
\end{split}
\end{equation*}
которому соответствует 2 "вещественных"\ соотношения --- равенство
нулю вещественной и мнимой части. Следовательно, среди параметров
(\ref{CD17}) только $8=|Q_v|+1$ независимых.

Покажем, что основой $D$ ортоскалярное представление $T$
определяется однозначно.

$$
\left[\begin{array}{c}
                                            A_1 \\
                                            A_2 \\
                                         \end{array}\right]=
\begin{array}{|cc|}
                                           \hline 0 & x_0 \\
                                           \hline x_1 & 0 \\
                                           0 & a_{22} \\
                                           0 & a_{32}\\
                                           \hline
\end{array}\ , \quad x_0>0.
$$

Из равенства длин столбцов следует
$x_0^2+a_{22}^2+a_{32}\overline{a_{32}}=x_1^2$, так что
положительное число $x_0$ по основе представления определяется
однозначно.
$$
\left[\begin{array}{c}
                                            B_1 \\
                                            B_2 \\
                                         \end{array}\right]=
\begin{array}{|cc|}
                                           \hline b_{01} & b_{02} \\
                                           \hline 0 & b_{12} \\
                                           b_{21} & 0 \\
                                           b_{31} & b_{32}\\
                                           \hline
\end{array}\ , \quad
\left[\begin{array}{c}
                                            C_1 \\
                                            C_2 \\
                                         \end{array}\right]=
\begin{array}{|cc|}
                                           \hline c_{01} & c_{02} \\
                                           \hline 0 & c_{12} \\
                                           c_{21} & 0 \\
                                           c_{31} & c_{32}\\
                                           \hline
\end{array}
$$

Унитарными преобразованиями строк элементы $b_{01}$ и $c_{01}$ можно
сделать вещественными неотрицательными:
$$
    b_{01}=y_{0}, \quad b_{02}=y_1e^{i\theta}.
$$
Из условий ортоскалярности вытекает, что
\begin{equation}
    \begin{split}
        b_{01}^2+b_{21}^2+|b_{31}|^2&=|b_{02}|^2+b_{12}^2+|b_{32}|^2,
        \quad \mbox{так что} \\
        b_{01}^2-|b_{02}|^2&=b_{21}^2+|b_{32}|^2-b_{21}^2-|b_{31}|^2=s.
        \\
        b_{01}b_{02}+\overline{b_{31}}b_{32}&=0,   \quad \mbox{так что} \\
        b_{01}b_{02}=b_{01}|b_{02}|e^{i\theta}&=\sin\phi_3\cos\psi_4\sin\psi_3\cos\phi_1\cos\psi_3e^{i(\theta_1-\theta_3)}.\\
    \end{split}
\end{equation}

Поэтому $\theta=\theta_1-\theta_3$ и
\begin{equation}
\begin{split}
    b_{01}^2\cdot(-|b_{02}|^2)=-\sin^2\phi_3\cos^2\psi_4\sin^2\psi_3\cos^2\phi_1\cos^2\psi_3=-t
\end{split}
\end{equation}
Стало быть числа $b_{01}^2$ и $-|b_{02}|^2$ есть корни квадратного
уравнения
$$
    z^2-s\cdot z-t=0,\quad t>0.
$$
Уравнение имеет два вещественных корня разного знака, таким образом
$|b_{01}|$, $|b_{02}|$ и $\theta$ определяются основой предствления
однозначно.

Аналогично по основе представления однозначно определяются и числа
$c_{01}$, $c_{02}$. Стало быть в размерности $\delta_{\widetilde
E_6}$ неразложимые ортоскалярные представления зависят не более чем
от $|Q_v|+1$, а неразложымые ортоскалярные представления общего
положения в точности от $8=|Q_v|+1$ параметров.

\medskip

d) Пусть $G=\widetilde E_7$; соответствующий разделённый колчан
имеет вид

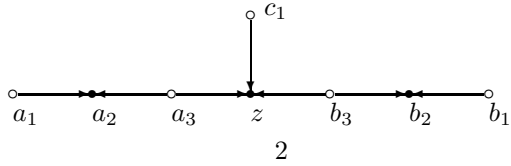
\begin{center}
\begin{picture}(200,40)

\put(0,0){\circle{3}} \put(30,0){\circle*{3}} \put(60,0){\circle{3}}
\put(90,0){\circle*{3}} \put(120,0){\circle{3}}
\put(150,0){\circle*{3}} \put(180,0){\circle{3}}
\put(90,30){\circle{3}}

\put(2,0){\vector(1,0){27}} \put(58,0){\vector(-1,0){27}}
\put(62,0){\vector(1,0){27}} \put(90,28){\vector(0,-1){27}}
\put(118,0){\vector(-1,0){27}} \put(122,0){\vector(1,0){27}}
\put(178,0){\vector(-1,0){27}}

\put(0,-10){$a_1$} \put(30,-10){$a_2$} \put(60,-10){$a_3$}
\put(90,-10){$z$} \put(120,-10){$b_3$}  \put(150,-10){$b_2$}
\put(180,-10){$b_1$} \put(95,30){$c_1$}
\end{picture}
\end{center}
\medskip
\medskip
\medskip

$
\begin{array}{ccccccccc}
                                     &   &   &   & 2 &   &  & \\
         \mbox{Как извесно,}\quad \delta_{\widetilde E_7}=   & 1 & 2 & 3 & 4 & 3 & 2& 1
                                   \end{array}
$

Пусть $T$ --- неразложимое ортоскалярное представление колчана $Q$ в
размерности $\delta_{\widetilde E_7}$,
\begin{gather*}
  T=
\left[\begin{array}{ccccc}
           A_{11} & A_{12} & 0 & 0 & 0 \\
           0 & A_{22} & A_{23} & A_{24} & 0 \\
           0 & 0 & 0 & A_{34} & A_{35} \\
         \end{array}\right]\\
          =
        \left[\begin{array}{c|c|c|c|c}
           \begin{array}{c}
             a_{11} \\
             a_{21}
           \end{array}
            & \begin{array}{ccc}
                a_{12} & a_{13} & a_{14} \\
                a_{21} & a_{22} & a_{24}
              \end{array}
             & 0 & 0 & 0 \\
                       \hline
          0 & \begin{array}{ccc}
                 a_{32} & a_{33} & a_{34} \\
                 a_{42} & a_{43} & a_{44} \\
                 a_{52} & a_{53} & a_{54} \\
                 a_{62} & a_{63} & a_{64}
               \end{array}
            & \begin{array}{cc}
                a_{35} & a_{36} \\
                a_{45} & a_{46} \\
                a_{55} & a_{56} \\
                a_{65} & a_{66}
              \end{array}
             & \begin{array}{ccc}
                 a_{37} & a_{38} & a_{39} \\
                 a_{47} & a_{48} & a_{49} \\
                 a_{57} & a_{58} & a_{59} \\
                 a_{67} & a_{68} & a_{69}
               \end{array} & 0 \\
                       \hline
           0 & 0 & 0 & \begin{array}{ccc}
                 a_{77} & a_{78} & a_{79} \\
                 a_{87} & a_{88} & a_{89}
               \end{array} &  \begin{array}{c}
                a_{7,10} \\
                a_{8,10}
              \end{array}\\
         \end{array}\right] \ ,
\end{gather*}

\begin{equation*}
\begin{split}
    A_{11}&=T_{a_2,a_1}; \ A_{12}=T_{a_2,a_3}; \ A_{22}=T_{z,a_3}; \
    A_{23}=T_{z,c_1}; \\
    A_{24}&=T_{z,b_3}; \ \ \; A_{34}=T_{b_2,b_3}; \ A_{35}=T_{b_2,b_1}.
\end{split}
\end{equation*}

Допустимыми унитарными преобразованиями строк и столбцов некоторые
из элементов матрицы $T$ преобразуем в нули, некоторые из елементов
станут равными нулю в силу ортогональности строк и столбцов внутри
полос.

Символ $\overline{0}_k$ на каком-либо место матрицы $T$ будет
означать, что на этом месте может буть получен ноль на $k$-м шаге
унитарными преобразованиями столбцов вертикальной полосы, $0|_k$ ---
строк горизонтальной полосы, $\overrightarrow{{0}}_k$ --- ноль
получен в силу ортогональности столбцов вертикальной полосы,
${0\!\!\downarrow_k}$ --- в силу ортогональности строк. При этом
делая нули на $k$-м шаге, мы не "портим"\ нули, полученные ранее.
Умножением строк и столбцов на числа вида $e^{i\phi}$ (унитарные
преобразования) некоторые из элементов можно сделать вещественными
неотрицательными, и мы на это укажем прямо в матрице ($a_{ij}\geq0$
либо $a_{ij}\leq0$). Мы указываем на строгое сравнение ($a_{ij}>0$
либо $a_{ij}<0$), если из равенства элемента $a_{ij}$ нулю следует
разложимость представления. В результате матрицу $T$ можно привести
к виду

\begin{gather*}
\left[
      \begin{array}{ccc}
        A_{11} & A_{12} & 0 \\
        0 & A_{22} & A_{23} \\
        0 & 0 & 0 \\
      \end{array}
    \right]= \\
\left[\begin{array}{c|c|c}
           \begin{array}{c}
             a_{11}\geq 0 \\
             a_{21}\geq 0 \\
           \end{array}
            & \begin{array}{ccc}
                0|_{10} & a_{13}\geq0 & a_{14}\leq0 \\
                a_{22}>0 & a_{23}\leq0 & \overrightarrow{{0}}_{11}
              \end{array}
             &\begin{array}{cc}
                   0\quad & \quad 0 \\
                   0 \quad & \quad 0 \\
                 \end{array} \\
                       \hline
          \begin{array}{c}
                   0  \\ 0  \\0  \\0  \\
                 \end{array} & \begin{array}{ccc}
                 a_{32}>0 & \overline{0}_3 & \overline{0}_3 \\
                 a_{42}>0 & a_{43}\geq0 & \overline{0}_7 \\
                 0|_4 & a_{53}\geq0 & a_{54}\geq0 \\
                 0|_4 & 0|_8 & a_{64}\geq0
               \end{array}
            & \begin{array}{cc}
                a_{35}>0 & \overline{0}_3 \\
                a_{45}<0 & \overrightarrow{{0}}_6 \\
               0\!\!\downarrow_5 & a_{56}\geq0 \\
                0\!\!\downarrow_5 & a_{66}
              \end{array}  \\
            \hline
           \begin{array}{cc}
                   0 \\
                   0 \\
                 \end{array} & \begin{array}{ccc}
                   0\quad \quad \quad & 0 & \quad \  \quad 0 \\
                   0\quad \quad \quad & 0 & \quad \ \quad 0 \\
                 \end{array} & \begin{array}{cc}
                   0\quad &\quad 0 \\
                   0\quad &\quad 0 \\
                 \end{array}
         \end{array}\right]
\end{gather*}

$$
\left[
  \begin{array}{cc}
    0 & 0 \\
    A_{24} & 0 \\
    A_{34} & A_{35} \\
  \end{array}
\right]=\left[
\begin{array}{c|c}
              \begin{array}{ccc}
                   0 \quad  \; & \quad\;  0 &\quad \quad  0 \\
                   0 \quad \; & \quad\; 0 &\quad \quad  0 \\
                 \end{array} & \begin{array}{c}
                   0  \\
                   0  \\
                 \end{array} \\
                       \hline
             \begin{array}{ccc}
                 \overrightarrow{{0}}_3 & \overrightarrow{{0}}_3 & a_{39} \\
                 a_{47}>0 & \overline{0}_7 & 0|_2 \\
                 a_{57}\leq0 & a_{58}\geq0 & 0|_2 \\
                 0\!\!\downarrow_9 & a_{68} & 0|_2
               \end{array} &    \begin{array}{c}
                   0  \\
                   0  \\
                   0  \\
                   0
                 \end{array} \\
                       \hline
            \begin{array}{ccc}
                 0|_{10} & a_{78}\geq0 & \overline{0}_1 \\
                 a_{87}\geq0 & a_{88}\geq0 & \overline{0}_1 \\
               \end{array} &  \begin{array}{c}
                a_{7,10}\geq0 \\
                a_{8,10}\leq0 \\
              \end{array}\\
         \end{array}\right]
$$
Нормируем представление $T$ (будем считать, что $\chi_z=1$). Тогда
матрицу $D=\left[
             \begin{array}{ccc}
               A_{22} & A_{23} & A_{24} \\
             \end{array}
           \right]$ (будем называеть ее основой представления) можно
параметризировать следующим образом:
\begin{gather*}
        A_{22}=\left[
        \begin{array}{ccc}
          \cos\phi_1\cos\psi_1 & 0 & 0   \\
          \sin\phi_1\cos\psi_2 &  \cos\phi_2\sin\psi_2 & 0  \\
          0 & \sin\phi_2\cos\psi_4 & \cos\phi_3\cos\psi_3\sin\psi_4\\
          0 & 0 & \sin\phi_4\sin\psi_5 \\
        \end{array}
      \right] \\
        A_{23}=
    \left[
        \begin{array}{cc}
           \sin\phi_1\cos\psi_1 & 0  \\
           -\cos\phi_1\cos\psi_2 & 0  \\
           0 & \sin\phi_3\cos\psi_3\sin\psi_4\\
           0 & \cos\phi_4\sin\psi_5e^{i\theta_1}\\
        \end{array}
      \right]\\
        A_{24}=
    \left[
        \begin{array}{ccc}
          0 & 0 & \sin\psi_1 \\
          \sin\phi_2\sin\psi_2 & 0 & 0 \\
          -\cos\phi_2\cos\psi_4 & \sin\psi_3\sin\phi_4 & 0 \\
          0 & \cos\psi_5e^{i\theta_2} & 0 \\
        \end{array}
    \right]
\end{gather*}

Строки матрицы $D$ ортонормированы, если дополнительно выполняется
соотношение (означающее ортогональность 3-й и 4-й строк):
\begin{equation} \label{CD19}
\begin{split}
    \cos\phi_3\sin\phi_4\cos\psi_3\sin\psi_4\sin\psi_5&+\sin\phi_3\cos\phi_4\cos\psi_3\sin\psi_4\sin\psi_5e^{i\theta_1}\\
&+\sin\psi_3\sin\psi_4\cos\psi_5e^{i\theta_2}=0.
\end{split}
\end{equation}
Это соотношение эквивалентно двум вещественным, так что из 11
параметров $\phi_1-\phi_4,\ \psi_1-\psi_5,\ \theta_1,\ \theta_2$
независимыми являются $9=|Q_v|+1$.

Легко видеть, что по основе $D$ представления $T$ матрицы $A_{11}$,
$A_{12}$, $A_{34}$, $A_{35}$ (ранее частично приведенные) находятся
однозначно.

Так из ортоскалярности столбцов матрицы $\left[
                                           \begin{array}{c}
                                             A_{24} \\
                                             A_{34} \\
                                           \end{array}
                                         \right]$ (т.е. из равенства
$A_{24}^*A_{24}+ A_{34}^*A_{34}=\chi_{b_3}I_{b_3}$) следует, что
\begin{equation*}
\begin{split}
    &a_{87}^2+a_{47}^2+a_{54}^2=a_{39}^2\\
    &a_{87}a_{88}+a_{54}a_{58}=0\\
    &a_{78}^2+a_{88}^2+a_{58}^2+a_{68}\overline{a_{68}}=a_{39}^2.
\end{split}
\end{equation*}
Из этих равенств последовательно и однозначно находятся элементы
$a_{87}$, $a_{88}$, $a_{78}$. Из ортоскалярности строк матрицы
 $\left[
                                           \begin{array}{cc}
                                             A_{34} &
                                             A_{35} \\
                                           \end{array}
                                         \right]$
(т.е. из равенства $A_{34}A_{34}^*+
A_{35}A_{35}^*=\chi_{b_2}I_{b_2}$) следует, что
\begin{equation*}
\begin{split}
    &a_{78}^2+a_{7,10}^2=a_{87}^2+a_{88}^2+a_{8,10}^2\\
    &a_{78}a_{88}+a_{7,10}a_{8,10}=0,
\end{split}
\end{equation*}
откуда однозначно находятся элементы $a_{7,10}$ и $a_{8,10}$.

Параметризуем столбцы матрицы $\left[
                                           \begin{array}{c}
                                             A_{12} \\
                                             A_{22} \\
                                           \end{array}
                                         \right]$ следующим образом:
$$
    \left[
                                           \begin{array}{c}
                                             A_{12} \\
                                             A_{22} \\
                                           \end{array}
                              \right]=
    \left[
      \begin{array}{ccc}
        0 & x\sin\alpha_2\sin\beta_2 & -x\sin\beta_3 \\
        x\sin\beta_1 & -x\cos\alpha_2\sin\beta_2 & 0 \\
        \hline
        x\cos\alpha_1\cos\beta_1 & 0 & 0 \\
        x\sin\alpha_1\cos\beta_1 & x\sin\alpha_3\cos\beta_2 & 0 \\
        0 & x\cos\alpha_3\cos\beta_2 & x\sin\alpha_4\cos\beta_3 \\
        0 & 0 & x\cos\alpha_4\cos\beta_3
      \end{array}
    \right]\ ,
$$
здесь $x>0$.

Покажем, что элементы матрицы $A_{12}$ элементами матрицы $A_{22}$
определяются однозначно. Так как $x\cos\alpha_1\cos\beta_1=a_{32}$,
$x\sin\alpha_1\cos\beta_1=a_{42}$, то
$x^2\cos^2\beta_1=a_{32}^2+a_{42}^2$. Аналогично,
$x^2\cos^2\beta_2=a_{43}^2+a_{53}^2$ и
$x^2\cos^2\beta_3=a_{54}^2+a_{64}^2$. Поэтому
\begin{equation*}
\begin{split}
    x\sin\beta_1&=\sqrt{x^2-(a_{32}^2+a_{42}^2)} \\
    x\sin\beta_2&=\sqrt{x^2-(a_{43}^2+a_{53}^2)} \\
    x\sin\beta_3&=\sqrt{x^2-(a_{54}^2+a_{64}^2)} \\
\end{split}
\end{equation*}
Из ортогональности столбцов матрицы $\left[
                                           \begin{array}{c}
                                             A_{12} \\
                                             A_{22} \\
                                           \end{array}
                                         \right]$ следует
\begin{equation*}
\begin{split}
    -x^2\sin\beta_1\sin\beta_2\cos\alpha_2+a_{42}a_{43}&=0 \\
    -x^2\sin\beta_2\sin\beta_3\sin\alpha_2+a_{53}a_{54}&=0
\end{split}
\end{equation*}
Поэтому
\begin{equation*}
\begin{split}
    \cos\alpha_2&=\frac{a_{42}a_{43}}{\sqrt{\big[x^2-(a_{32}^2+a_{42}^2)\big]\cdot\big[x^2-(a_{43}^2+a_{53}^2)\big]}}
    \\
    \sin\alpha_2&=\frac{a_{53}a_{54}}{\sqrt{\big[x^2-(a_{54}^2+a_{64}^2)\big]\cdot\big[x^2-(a_{43}^2+a_{53}^2)\big]}}
\end{split}
\end{equation*}

Значение $x=\chi_{a_3}$ определяется по основе представления
однозначно: если $A_{22}$ унитарно эквивалентна матрице
$\textrm{diag}\{r_1,r_2,r_3\}$, $r_i\geq0$, то
$x=\max\limits_{i}r_i$ (это легко получить, диагонализируя матрицы
$A_{12}$, $A_{22}$ унитарными преобразованиями, при этом существенно
используется ортоскалярность представления). Элементы $a_{11}$,
$a_{21}$ однозначно находятся по матрице $A_{12}$ (аналогично тому,
как находятся єлементы $a_{7,10}$, $a_{8,10}$ по матрице $A_{34}$).

\medskip

e) Пусть $G=\widetilde E_8$; соответствующий разделенный колчан
имеет вид

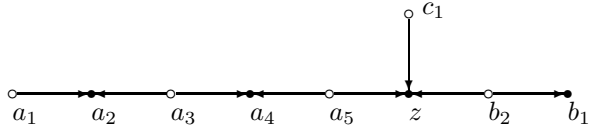
\begin{center}
\begin{picture}(220,40)

\put(0,0){\circle{3}} \put(30,0){\circle*{3}} \put(60,0){\circle{3}}
\put(90,0){\circle*{3}} \put(120,0){\circle{3}}
\put(150,0){\circle*{3}} \put(180,0){\circle{3}}
\put(150,30){\circle{3}} \put(210,0){\circle*{3}}

\put(2,0){\vector(1,0){27}} \put(58,0){\vector(-1,0){27}}
\put(62,0){\vector(1,0){27}} \put(150,28){\vector(0,-1){27}}
\put(118,0){\vector(-1,0){27}} \put(122,0){\vector(1,0){27}}
\put(178,0){\vector(-1,0){27}} \put(182,0){\vector(1,0){27}}

\put(0,-10){$a_1$} \put(30,-10){$a_2$} \put(60,-10){$a_3$}
\put(90,-10){$a_4$} \put(120,-10){$a_5$}  \put(150,-10){$z$}
\put(180,-10){$b_2$}  \put(210,-10){$b_1$} \put(155,30){$c_1$}
\end{picture}
\end{center}

\medskip
\medskip

Как извесно,
$$
\begin{array}{cccccccccc}
                                     &   &   &   &   &   & 3 &   &\\
          \delta_{\widetilde E_8}=   & 1 & 2 & 3 & 4 & 5 & 6 & 4 & 2
                                   \end{array}
$$

Пусть $T$ --- неразложимое ортоскалярное представление колчана $Q$ в
размерности $\delta_{\widetilde E_8}$,
$$
    T=\left[\begin{array}{ccccc}
        A_{11} & A_{12} & 0 & 0 & 0 \\
        0 & A_{22} & A_{23} & 0 & 0 \\
        0 & 0 & A_{33} & A_{34} & A_{35} \\
        0 & 0 & 0 & 0 & A_{45} \\
      \end{array}
      \right],
$$
где
\begin{equation*}
\begin{split}
    A_{11}&=T_{a_2,a_1}; \ A_{12}=T_{a_2,a_3};\ A_{22}=T_{a_4,a_3}; \ A_{23}=T_{a_4,a_5};\\
    A_{33}&=T_{z,a_5}; \ A_{34}=T_{z,c_1}; \ A_{35}=T_{z,b_2}; \
    A_{45}=T_{b_1,b_2}.
\end{split}
\end{equation*}

Допустимыми унитарными преобразованиями строк и столбцов некоторые
мз елементов матрицы $T$ преобразуем в нули, некоторые из елементов
станут равными нулю в силу ортогональность строк и стоблцов внутри
полос. Следуя предыдущим договорённостям об обозначениях, матрицу
$T$ можно привести к виду, при котором

\begin{gather*}
\left[
  \begin{array}{ccc}
    A_{11} & A_{12} \\
    0 & A_{22}  \\
  \end{array}
\right]=  \left[
  \begin{array}{c|ccc}
    a_{11}>0 & a_{12}>0 &  \overrightarrow{{0}}_{26}  & \overrightarrow{{0}}_{26}\\
    0\!\!\downarrow_{27} & O|_{25} & a_{23}>0 & a_{24}<0\\
    \hline
    0 & a_{32}>0 & \overline{0}_{21} & \overline{0}_{21}\\
    0 & 0\!\!\downarrow & a_{43}>0 &  \overline{0}_{23}  \\
    0 & 0\!\!\downarrow & a_{53}>0 & a_{54}>0 \\
    0 & 0\!\!\downarrow & 0\!\!\downarrow_{24} & a_{64}>0 \\
  \end{array}
\right]
\\
\left[
  \begin{array}{c}
    A_{23} \\
    A_{33} \\
  \end{array}
\right]=  \left[
  \begin{array}{ccccc}
     a_{35}>0 & \overrightarrow{{0}}_{17} & \overrightarrow{{0}}_{17} & \overrightarrow{{0}}_{17} & \overrightarrow{{0}}_{17} \\
     0|_{16} & a_{46}>0 & a_{47}<0 & \overrightarrow{{0}}_{19} & \overrightarrow{{0}}_{19} \\
     0|_{16} & 0|_{18} & a_{57}>0 & a_{58}<0 & \overrightarrow{{0}}_{21} \\
     0|_{16} & 0|_{18} & 0|_{20} & a_{59}>0 & a_{69}<0 \\
     \hline
     a_{75}>0 & \overrightarrow{{0}}_2 & \overrightarrow{{0}}_2 & \overrightarrow{{0}}_2 & \overrightarrow{{0}}_2 \\
     0|_1 & a_{86}\geq0 & \overline{0}_6 & \overline{0}_6 & \overline{0}_6 \\
     0|_1 & a_{96}\geq0 & a_{97}\geq0 & \overline{0}_7 & \overline{0}_7 \\
     0|_1 & 0|_7 & a_{10,7}\geq0 & a_{10,8}\geq0 & \overline{0}_{11} \\
     0|_1 & 0|_7 & 0|_{10} & a_{11,8}\geq0 & a_{11,9}\geq0 \\
     0|_1 & 0|_7 & 0|_{10} & a_{12,9}\geq0 & a_{12,9}\geq0
  \end{array}
\right] \end{gather*}

\begin{gather*}
\left[
  \begin{array}{cc}
    A_{34} & A_{35}\\
    0 & A_{45}\\
  \end{array}
\right]=\\
\left[
  \begin{array}{ccc|ccccc}
         a_{7,10}>0 & \overline{0}_2 & \overline{0}_2 & a_{7,13}>0 & \overline{0}_2 & \overline{0}_2 & \overline{0}_2 \\
         a_{8,10}>0 & \overrightarrow{{0}}_6 &
         \overrightarrow{{0}}_6
         &          a_{8,13}<0 & a_{8,14}>0 & \overrightarrow{{0}}_6 & \overrightarrow{{0}}_6 \\
         0|_3 & \!\!\!\! a_{9,11}>0 & \overline{0}_7 &           0\!\!\downarrow_4 & a_{9,14}<0 & \overrightarrow{{0}}_9 & \overrightarrow{{0}}_9 \\
         0|_3 & \!\!\!\! a_{10,11}<0 & \overline{0}_{13} &           0\!\!\downarrow_4 & 0\!\!\downarrow_8 & a_{10,15}>0 & \overline{0}_{14} \\
         0|_3 & 0\!\!\downarrow_{14} &\!\!\!\! a_{11,12}\geq0 &           0\!\!\downarrow_4 & 0\!\!\downarrow_8 & a_{11,15}<0 & a_{11,16}\geq0 \\
         0|_3 & 0\!\!\downarrow_{14} &\!\!\!\! a_{12,12} &           0\!\!\downarrow_4 & 0\!\!\downarrow_8 & 0\!\!\downarrow_{15} &
          a_{12,16} \\
          \hline
          0 & 0 &0 & a_{13,13}>0 & a_{13,14}>0 & 0_5 & 0_5 \\
          0 & 0 &0  &          0\!\!\downarrow_5 & 0\!\!\downarrow_5 & a_{14,15}>0 & a_{14,16}>0 \\
        \end{array}
\right]
\end{gather*}

Основу представления $D=\left[
                          \begin{array}{ccc}
                            A_{33} & A_{34} & A_{35} \\
                          \end{array}
                        \right]$ можно параметризовать следующим
                        образом
\begin{equation*}
\begin{split}
    a_{75}&=\sin\psi_1, \ a_{7,10}=\cos\phi_1\cos\psi_1, \
    a_{7,13}=\sin\phi_1\cos\psi_1 \\
    a_{86}&=\sin\phi_2\sin\psi_2, \ a_{8,10}=\sin\phi_1\cos\psi_2,\
    a_{8,13}=-\cos\phi_1\cos\psi_2,\\ &a_{8,14}=\cos\phi_2\sin\psi_2\\
    a_{96}&=\cos\phi_2\sin\psi_3, \ a_{97}=\sin\phi_3\cos\psi_3, \
    a_{9,11}=\cos\phi_3\cos\psi_3; \\
    &a_{9,14}=-\sin\phi_2\sin\psi_3
\end{split}
\end{equation*}
\begin{equation*}
\begin{split}
  a_{10,7}&=\cos\phi_3\sin\psi_4; \ a_{10,8}=\sin\phi_4\cos\psi_4;
    \ a_{10,11}=-\sin\phi_3\sin\psi_4; \\
    &a_{10,15}=-\cos\phi_4\cos\psi_4 \\
    a_{11,8}&=\cos\phi_4\cos\psi_4\cos\psi_5; \
    a_{11,9}=\cos\phi_5\sin\psi_{5}; \
    a_{11,12}=\sin\phi_5\sin\psi_5; \\
    &a_{11,15}=\sin\phi_4\cos\psi_4\cos\psi_5;\
    a_{11,16}=\sin\psi_4\cos\psi_5 \\
    a_{12,9}&=\cos\phi_6\cos\psi_6;\
    a_{12,12}=\sin\phi_6\cos\psi_6e^{i\theta_1}; \
    a_{12,16}=-\sin\psi_6e^{i\theta_2}.
\end{split}
\end{equation*}

Строки матрицы $D$ ортонормированы, если дополнительно выполняется
соотношение (означающее ортогональность 5-й и 6-й строк):
\begin{equation*}
\begin{split}
    \cos\phi_5\cos\phi_6\sin\psi_5\cos\psi_6&+\sin\phi_5\sin\phi_6\sin\psi_5\cos\psi_6e^{i\theta_1}\\
    &-\sin\psi_4\cos\psi_5\sin\psi_6e^{i\theta_2}=0.
\end{split}
\end{equation*}

Это соотношение эквивалентно двум вещественным.

Кроме того из ортоскалярности представления следует равенство длин
столбцов (как векторов в унитарном пространстве) матрицы $A_{34}$,
поэтому имеем еще 2 соотношения:
\begin{equation*}
\begin{split}
    \cos^2\phi_1\cos^2\psi_1+\sin^2\phi_1\cos^2\psi_2&=\cos^2\phi_3\cos^2\psi_3+\sin^2\phi_3\sin^2\psi_4\\
\mbox{и} \quad
    \cos^2\phi_1\cos^2\psi_1+\sin^2\phi_1\cos^2\psi_2&=\sin^2\phi_5\sin^2\psi_5+\sin^2\phi_6\cos^2\psi_6.
\end{split}
\end{equation*}
Таким образом из 14 параметров $\phi_1-\phi_6, \ \psi_1-\psi_6,\
\theta_1,\ \theta_2$ независимыми являются $10=|Q_v|+1$.

Также, как и в предыдущих случаях, легко показать, что по основе $D$
представления $T$ матрицы $A_{11},\ A_{12},\ A_{22},\ A_{23},\
A_{45}$ (ранее частично приведенные) находятся однозначно.
\endproof

\begin{rem}
    Неразложимые представления расширенных графом Дынкина в
    категории линейных пространств изучены в \cite{Naz73},\cite{DonFre73}.

    Как и для категории линейных пространств, размерности
    неразложимых ортоскалярных представлений этих графов
    не ограничены в совокупности. В размерности, являющейся
    действительным корнем графа, при фиксированном характере,
    неразложимое представление единственно. Для размерностей,
    являющихся мнимыми корнями графа, неразложимые ортоскалярные
    представления, в отличие от линейно-алгебраического случая,
    существуют только в размерности совпадающей с минимальным
    положительным мнимым корнем, но также зависят, при фиксированном
    характере, от одного  (комплексного) параметра.
\end{rem}




\end{document}